\documentstyle[11pt]{article}

\title{BOUNDS ON THE NODAL STATUSES OF SOME TRANSFINITE GRAPHS}
\author{A. H. Zemanian}
\date{}

\begin{document}
\maketitle
\baselineskip21pt

{\ Abstract --- The bounds on the statuses of the nodes in a finite graph
established by Entringer, Jackson, and Snyder 
are extended herein so that they apply to the statuses  
of the nodes in transfinite graphs of a certain kind.
\\

} 

\section{Introduction}

The purpose of this note is to extend the known bounds on the
statuses of the nodes of a finite graph to nodes in transfinite graphs 
of a certain kind.  That known result was established by
Entringer, Jackson, and Snyder \cite{ejs}.\footnote{See also
\cite[pages 43-44]{b-h} for an exposition of that result.}  
It states that the status $s(x)$
of any node $x$ in a finite connected graph $G$ having $p$ nodes 
and $q$ branches satisfies the inequalities
\begin{equation}
p-1\;\leq\;s(x)\;\leq\;(p-1)(p+2)/2\,-\,q  \label{1.1}
\end{equation}
and that these bounds can be achieved for each $q$ such that 
$(p-1)\;\leq q\;\leq\;p(p-1)/2$.  A modification of this result holds for 
transfinite graphs satisfying certain conditions.

\section{Some Preliminary Definitions and Known Results}

We shall use some definitions and symbolism appearing in \cite{b8}.
Also, we restrict our attention to transfinite graphs 
${\cal G}^{\mu}$ of rank $\mu$,
where for the sake of some simplicity we restrict $\mu$ to the 
positive natural 
numbers.\footnote{When $\mu=0$, ${\cal G}^{\mu}$ is a conventional
graph.}  The rather complicated recursive 
definitions of such $\mu$-graphs appear in Section 2.4 of 
\cite{b8}.  All our arguments extend readily to graphs of 
higher ranks, that is, the transfinite-ordinal ranks. Since the 
distance from any nonmaximal node $z$ 
is the same as the distance from any node of higher rank containing $z$,
we can restrict our attention to 
the maximal nodes in ${\cal G}^{\mu}$, that is, 
to the nodes that are not contained in any node of higher rank.  
This is understood henceforth.

Transfinite nodes
are defined in terms of tips (i.e., graphical extremities),
which in turn are equivalence classes of one-ended paths, 
as stated in \cite[page 11]{b8}.  
The $\mu$-nodes in ${\cal G}^{\mu}$ 
are the nodes of highest rank in ${\cal G}^{\mu}$.
A $\mu$-node is said 
to be pristine if it does not contain a node of lower rank;  
We will assume that all the $\mu$-nodes are pristine.
Also,
a node $x^{\rho}$ of any rank $\rho$ $(\rho\leq \mu)$ is called 
a nonsingleton if it contains at least two elements 
(either two $(\mu-1)$-tips of a $(\mu-1)$-tip and a node of 
rank lower than $\rho$).

Furthermore, two branches are said to be $\rho$-connected
if there is a path of rank $\rho$ or less that terminates at 
nodes of those branches.  Actually, such path-connectedness need 
not exist between all pairs of branches.\footnote{A more general 
concept of connectedness is based on walks.  Such walk-connectedness
always exists between branches.  Our results extend to this case 
\cite[page 67]{b8}, as is indicated at the end of this paper.}
To insure that such path-connectedness does exist, we impose
the following Condition A \cite[page 25]{b8}.  We say that two tips 
are nondisconnectable if their representative paths meet infinitely often
\cite[page 25]{b8}.  Also, a node is said to embrace
a tip if that tip is part of that node (see \cite[page 12]{b8} 
for the precise definition).  

{\bf Condition A.} If two tips are nondisconnectable, then either
they are contained in the same node or at least one of them 
is the sole member of a maximal node.  

Under this condition, for any two nonsingleton nodes there will be a path
that terminates at them \cite[Lemma 4.3-2]{b8}, and moreover
such path connectedness is a transitive binary relation 
for the nonsingleton nodes of ${\cal G}^{\mu}$;  in fact, it is an equivalence
relation \cite[Theorem 3.1-4]{b8}.

Throughout this work we assume that ${\cal G}^{\mu}$ is $\mu$-connected
in the sense that every pair of branches are $\rho$-connected for some rank 
$\rho\leq \mu$ depending on the choice of those branches.
As a result, the set of branches in ${\cal G}^{\mu}$ is partitioned
into subsets according to $(\mu-1)$-connectedness,
and the subgraph of rank $\mu-1$ induced by such a subset is called a 
$(\mu-1)$-section \cite[page 23]{b8}. Because we are assuming
that all the $\mu$-nodes are pristine, it follows that every 
node of rank less than $\mu$ is contained within some
$(\mu-1)$-section.  Thus, the $(\mu-1)$-sections also partition
the set of nodes of ranks less than $\mu$.
We say the a $\mu$-node $x^{\mu}$ is incident to a $(\mu-1)$-section
${\cal S}^{\mu-1}$
if $x^{\mu}$ contains a $(\mu-1)$-tip whose representative paths
lie within ${\cal S}^{\mu-1}$.  Moreover, if a $\mu$-node is incident to 
two or more $(\mu-1)$-sections, it serves as a connection
between them.

Lengths of paths and the distances between nodes are defined in 
\cite[Sections 4.2 and 4.4]{b8}. By virtue of Condition A, 
the $\mu$-connectedness of ${\cal G}^{\mu}$,
and the assumption that all the $\mu$-sections are pristine,
we have the following results as a consequence of \cite[Lemma 4.7-4]{b8}.
The length of any path $P$ within a $(\mu-1)$-section ${\cal S}^{\mu-1}$
that is incident to a $\rho$-node $x^{\rho}$
($\rho<\mu$) in ${\cal S}^{\mu-1}$
and reaches a $\mu$-node $x^{\mu}$ incident to 
${\cal S}^{\mu-1}$ is $\omega^{\mu}$.  This is because $P$ can reach $x^{\mu}$
only through a $(\mu-1)$-tip.  Consequently, 
we can define the $\mu$-length of any two-ended path $P$ that reaches or 
passes through at least one $\mu$-node as 
$\omega^{\mu}\cdot n$, where $n$ is the number of incidences that
$P$ makes with $\mu$-nodes; that is, when $P$ terminates at a 
$\mu$-node, there is one such incidence, and, when $P$
passes through a $\mu$-node from one $(\mu-1)$-section to another 
adjacent $(\mu-1)$-section, there are two such incidences. 
Furthermore, we define the $\mu$-distance between any two nodes $x$ and $y$ 
as the minimum of the $\mu$-lengths of all the 
paths that meet $x$ and $y$;  such a $\mu$-distance 
exists because those $\mu$-lengths 
are ordinals and any set of ordinals is well-ordered and therefore 
has a minimum.  It is a fact that under our assumptions
there will be a path terminating at $x$ and $y$ whose $\mu$-length is 
that $\mu$-distance;  such a path is called an $x,y$ geodesic.

\section{The $\mu$-Statuses of Nonsingleton Nodes}

Even though ${\cal G}^{\mu}$ is branchwise $\mu$-connected,
it can happen that there is no path between two nodes if 
at least one of them is a singleton, in which 
case no (path-based) distance will exist between them.
However, under Condition A, distances between nonsingleton 
nodes always exist.\footnote{See \cite[Section 3.1]{b8} for a 
discussion of this matter.}  For this reason, we shall restrict our 
definition of nodal statuses to the nonsingleton 
nodes.\footnote{One might motivate this restriction by noting that
a singleton node is a ``dead end'' in the sense that no path 
can pass through it, and so it does not contribute to the connectivity of 
${\cal G}^{\mu}$.}  

Another convention that we shall adopt in order to extend 
(\ref{1.1}) transfinitely is that the distance between any 
two nodes within a $(\mu-1)$-section is taken to be 0.
Thus, it is only a transition to or from a $\mu$-node that
contributes to the length of a geodetic path and thereby
to a distance.  Without this assumption, the status
of a node could be infinite.
Also, we shall henceforth assume that the $\mu$-connected 
$\mu$-graph ${\cal G}^{\mu}$
has only finitely many nonsingleton $\mu$-nodes 
and only finitely many $(\mu-1)$- sections.

To define the $\mu$-status of any nonsingleton node $x^{\rho}$
$(\rho\leq \mu)$, we first choose a single node $y_{m}^{\alpha_{m}}$
$(\alpha_{m}< \mu)$ for each $(\mu-1)$-section 
${\cal S}^{\mu-1}_{m}$, one such node for each $(\mu-1)$-section, 
and designate it as the representative node for 
${\cal S}^{\mu-1}_{m}$.  
Because we have taken the distances between nodes in a single 
$(\mu-1)$-section to be 0, we can take the distance from any node 
in ${\cal S}^{\mu-1}_{m}$ to 
be the same as the distance from the representative node 
$y_{m}^{\alpha_{m}}$ for ${\cal S}^{\mu-1}_{m}$.  
Then, we define the $\mu$-status
$s^{\mu}(x^{\rho})$ as the sum of the distances from $x^{\rho}$ to all the 
nonsingleton $\mu$-nodes plus the sum of the distances to the 
representative nodes of all the $(\mu-1)$-sections.  In symbols,
\begin{equation}
s^{\mu}(x^{\rho})\;=\;\sum_{k=1}^{K} d(x^{\rho},x_{k}^{\mu})\:+\:
\sum_{m=1}^{M} d(x^{\rho},y_{m}^{\alpha_{m}})  \label{3.1}
\end {equation}
Here, $k$ numbers the nonsingleton $\mu$-nodes, there being
$K$ of them, and $m$ numbers the $(\mu-1)$-sections, there being
$M$ of those.  By our assumptions, $K$ and $M$ are natural numbers.

\section{Conditions Imposed upon ${\cal G}^{\mu}$}

Let us now list all the conditions we have assumed for 
${\cal G}^{\mu}$.
\begin{description}
\item[4.1.]  The rank $\mu$ of the transfinite graph ${\cal G}^{\mu}$ is a 
positive natural number.
\item[4.2.]  ${\cal G}^{\mu}$ is $\mu$-connected (i.e., between every two 
branches there is path connecting them).
\item[4.3.]  Condition A holds.
\item[4.4.]  All $\mu$-nodes are pristine (i.e., none of them 
contains a node of lower rank).
\item[4.5.]  There are only finitely many nonsingleton $\mu$-nodes 
and only finitely many $(\mu-1)$-sections.
\item[4.6.]  The distance between any two nodes within the same
$(\mu-1)$-section is taken to be 0.
\end{description}

Also, bear in mind the following considerations. 

\noindent
{\bf Note 4.7.} We have restricted our attention to only the maximal
nodes because the distance from any nonmaximal node is 
the same as the distance 
from the maximal node containing it.

\noindent
{\bf Note 4.8.} Similarly, we have considered in our analysis 
only the nonsingleton nodes because the singleton nodes do not
contribute to the connectivity of ${\cal G}^{\mu}$.  On the other hand,
as a result of Assumptions 4.2 and 4.3, for any two 
nonsingleton nodes $x$ and $y$ there is a path terminating 
at $x$ and $y$ \cite[Lemma 4.3-2]{b8}.
Consequently, a distance is defined between $x$ and $y$.

In view of all this, the $\mu$-status $s(x^{\rho})$
for any maximal 
nonsingleton node $x^{\rho}$ $(\rho\leq\mu)$
is well-defined by (\ref{3.1}).

\section{The Replacement 0-graph}

In order to extend (\ref{1.1}) transfinitely, we replace ${\cal G}^{\mu}$
by a 0-graph in which the distances in ${\cal G}^{0}$ are the same as 
the $\mu$-distances in ${\cal G}^{\mu}$ except for a multiplicative factor 
$\omega^{\mu}$.  To do this, we adapt the replacement procedure given in
\cite[page 60]{b8}.  Remember that $\mu\geq 1$, that 
$k=1,2,\ldots,K$ numbers the nonsingleton $\mu$-nodes, and that
$m=1,2,\ldots,M$ numbers the $(\mu-1)$-sections.  Now,
we replace each nonsingleton $\mu$-node $x_{k}^{\mu}$ 
by a 0-node $x_{k}^{0}$.  Furthermore, having chosen a 
nonsingleton $\rho_{m}$-node $y_{m}^{\rho_{m}}$
$(0\leq \rho_{m}\leq \mu-1)$
within each $(\mu-1)$-section, we replace $y_{m}^{\rho_{m}}$
by a 0-node $y_{m}^{0}$.  (If $\rho_{m}=0$, no replacement is needed.)
Thus, for each $m$ we insert a branch between $y_{m}^{0}$ and
each $x_{k}^{0}$ corresponding to a $\mu$-node incident to the 
$(\mu-1)$-section ${\cal S}^{\mu-1}_{m}$ that contains $y_{m}^{\rho}$.
In this way, ${\cal G}^{\mu}$ is replaced by a finite 0-graph
consisting of the 0-nodes $x_{k}^{0}$ and $y_{m}^{0}$ $(k=1,2,\ldots,K;$
$m=1,2,\ldots, M)$ and of the said branches.  In particular,
each $(\mu-1)$-section ${\cal S}^{\mu-1}_{m}$ is replaced
by a star 0-graph with $y_{m}^{0}$ as its center 0-node and the 
$x_{k}^{0}$ corresponding to the $\mu$-nodes incident to
${\cal S}^{\mu-1}_{m}$ as its peripheral 0-nodes.

Now, to each path $P$ in ${\cal G}^{\mu}$ there corresponds a unique path $Q$
that passes through the 0-nodes $x_{k}^{0}$ corresponding 
to the $\mu$-nodes $x_{k}^{\mu}$.  If one terminal node $y_{m}^{\rho_{m}}$ 
$(\rho_{m} < \mu)$ of $P$
lies within a $(\mu-1)$-section, the path $Q$ terminates 
correspondingly at the 0-node 
$y_{m}^{0}$.  On the other hand, 
if $P$ terminates at a $\mu$-node $x_{k}^{\mu}$, then $Q$ 
terminates at $x_{k}^{0}$.  Furthermore, we define the $\mu$-length 
$|P|$ of $P$ as $\omega^{\mu}$ times the number $n$ of incidences 
that $P$ makes with the $\mu$-nodes, where 
we count an incidence once of $P$ terminates at a $\mu$-node 
and we count the incidence twice if $P$ passes through a 
$\mu$-node.  On the other hand, the number of branches in 
$Q$ is simply $n$. Thus, the length $|Q|$ of $Q$ is $n$.
As a result of all this, the lengths of $P$ and $Q$ are related as follows:
\begin{equation}
|P|\;=\;\omega^{\mu}\cdot |Q|  \label{5.1}
\end{equation}

\section{Bounds on the Statuses}

In conformity with the way we have defined the $\mu$-lengths of
paths in ${\cal G}^{\mu}$,
we can define the $\mu$-distance between any two nodes in 
${\cal G}^{\mu}$ as $\omega^{\mu}$ times the distance 
between the corresponding 0-nodes in ${\cal G}^{0}$, 
where any node within a $(\mu-1)$-section ${\cal S}^{\mu-1}_{m}$ is
represented by the node $y_{m}^{\rho_{m}}$ for 
${\cal S}^{\mu-1}_{m}$.  

Thus, the $\mu$-status of any 
node $x^{\rho}$ $(\rho\leq \mu)$ of ${\cal G}^{\mu}$, as defined by
(\ref{3.1}), is simply $\omega^{\mu}$ times the status of the 0-node in ${\cal G}^{0}$.

To lift the bounds (\ref{1.1}) to the transfinite case, we need merely determine the 
number $p$ of 0-nodes and the number $q$ of branches in ${\cal G}^{0}$.
Specifically, $p$ is the number 
of $\mu$-nodes in ${\cal G}^{\mu}$ plus the number of 
$(\mu-1)$-sections in ${\cal G}^{\mu}$.  With regard to $q$, note again that each 
$(\mu-1)$-section ${\cal S}^{\mu-1}_{m}$ has been replaced by a star graph with center
at $y_{m}^{0}$ and branches between $y_{m}^{0}$ and every one
of the $x_{k}^{0}$ corresponding to the $\mu$-nodes $x_{k}^{\mu}$
incident to ${\cal S}^{\mu-1}_{m}$.  Let $\delta_{m}$ be the number of such $x_{k}^{\mu}$;  
$\delta_{m}$ is the degree of $y_{m}^{0}$.  Then, 
\[ q\;=\;\sum_{m=1}^{M} \delta_{m}. \]
Altogether then, under the assumptions 4.1 to 4.6, we have the desired bounds 
on the status $s(x^{\rho})$ of any node $x^{\rho}$ in ${\cal G}^{\mu}$ as follows:
\begin{equation}
\omega^{\mu}\cdot (p-1)\;\leq\;s(x^{\rho})\;\leq\;\omega^{\mu}\cdot [(p-1)(p+2)/2\,-\,q]  \label{6.1}
\end{equation}
These bounds can be achieved for each $q$ such that $p-1\,\leq\,q\,\leq \,p(p-1)/2$.

\section{Two Final Comments}

{\bf 7.1.}  We have restricted our analysis to only the nonsingleton nodes in 
${\cal G}^{\mu}$.  Thus, the $\mu$-status $s^{\mu}(x)$ 
of any such node $x$ is the sum of 
the distances from $x$ to all (and only) the other nonsingleton 
$\mu$-nodes and the other representative nonsingleton nodes in 
the $(\mu-1)$-sections (one representative node to each
$(\mu-1)$-section.  This restriction can be relaxed somewhat 
by allowing distances from $x$ to include some of the singleton $\mu$-nodes.
Specifically, we can also allow a finite number of singleton 
$\mu$-nodes such that the $\mu$-distances to them from any nonsingleton node 
exists.  Such singleton nodes can occur;  see the set $\cal M$ defined in 
\cite[page 44]{b8}.

{\bf 7.2.}  We can also relax the restriction imposed by 
assumption 4.2 by taking definitions of entities in the transfinite graph 
to be based on walks rather than on paths.  As a result, two branches 
that are not connected by a transfinite 
path may be connected by a transfinite walk, as is explained in 
\cite[Chapter 5]{b8}.  Thus, upon assuming that the walk-based  
transfinite graph ${\cal G}^{\mu}$ is walk-connected, 
we can define walk-based 
distances between any two nodes of ${\cal G}^{\mu}$ 
\cite[Section 5.4]{b8}.  In this
case, we can discard assumption 4.3;  it is no longer needed.
Then, ${\cal G}^{\mu}$ can be related to a 
unique 0-graph ${\cal G}^{0}$ whereby the 
$\mu$-distance between nonsingleton walk-based nodes in 
${\cal G}^{\mu}$ is $\omega^{\mu}$ times the distance between
the corresponding 0-nodes in ${\cal G}^{0}$.
The procedure for doing this is much the same as that presented above.
As a result, we again have (\ref{3.1}) and 
(\ref{5.1}), and also (\ref{1.1})
replaced by (\ref{6.1}), where $s^{\mu}(x^{\rho})$ 
is now defined as the sum of the walk-based $\mu$-distances 
from $x^{\rho}$ to the walk-based nonsingleton nodes 
in ${\cal G}^{\mu}$.
Moreover, finitely many singleton $\mu$-nodes can also be 
allowed in this case, as is explained in the preceding 
paragraph 7.1.

\end{document}